\documentclass[12pt]{article}

\usepackage{geometry}                % See geometry.pdf to learn the layout options. There are lots.
\usepackage{graphicx}
\usepackage{amssymb}
\usepackage{epstopdf}
\usepackage{latexsym, amssymb, amsmath, amsthm, graphics}
\usepackage[mathscr]{eucal}
\usepackage{url}
\usepackage{eucal}
\title{The Amalgamated Product Structure of the Tame Automorphism Group in Dimension Three}
\author{David Wright}

\begin{document}
    \maketitle

 \input amssym.def
    \input amssym
    \theoremstyle{definition} 
    \newtheorem{rema}{Remark}[section]
    \newtheorem{questions}[rema]{Questions}
    \newtheorem{assertion}[rema]{Assertion}
    \newtheorem{num}[rema]{}
           \newtheorem{exam}[rema]{Example}
        \newtheorem{ques}[rema]{Question}
    \newtheorem{claim}[rema]{Claim}
    \theoremstyle{plain} 
    \newtheorem{propo}[rema]{Proposition}
    \newtheorem{theo}[rema]{Theorem} 
    \newtheorem{conj}[rema]{Conjecture}
    \newtheorem{quest}[rema]{Question} 
    \theoremstyle{definition}
    \newtheorem{defi}[rema]{Definition}
    \theoremstyle{plain}
    \newtheorem{lemma}[rema]{Lemma} 
    \newtheorem{corol}[rema]{Corollary}
    \newtheorem{rmk}[rema]{Remark}
    \newcommand{\del}{\triangledown}
    \newcommand{\nno}{\nonumber} 
    \newcommand{\lbar}{\big\vert}
    \newcommand{\mbar}{\mbox{\large $\vert$}} 
    \newcommand{\p}{\partial}
    \newcommand{\dps}{\displaystyle} 
    \newcommand{\bra}{\langle}
    \newcommand{\ket}{\rangle} 
    \newcommand{\kr}{\mbox{\rm Ker}\ }
    \newcommand{\res}{\mbox{\rm Res}} 
    \renewcommand{\hom}{\mbox{\rm Hom}}
    \newcommand{\pf}{{\it Proof:}\hspace{2ex}}
    \newcommand{\epf}{\hspace{2em}$\Box$}
    \newcommand{\epfv}{\hspace{1em}$\Box$\vspace{1em}}
    \newcommand{\nord}{\mbox{\scriptsize ${\circ\atop\circ}$}}
    \newcommand{\wt}{\mbox{\rm wt}\ } 
    \newcommand{\clr}{\mbox{\rm clr}\ }
    \newcommand{\ideg}{\mbox{\rm Ideg}\ } 
    \newcommand{\gC}{{\mathfrak g}_{\mathbb C}}
    \newcommand{\hatC}{\widehat {\mathbb C}} 
    \newcommand{\bZ}{{\mathbb Z}} 
    \newcommand{\bQ}{{\mathbb Q}}
    \newcommand{\bR}{{\mathbb R}} 
    \newcommand{\bN}{{\mathbb N}}
    \newcommand{\bT}{{\mathbb T}} 
    \newcommand{\fg}{{\mathfrak g}}
    \newcommand{\fgC}{{\mathfrak g}_{\bC}} 
    \newcommand{\cD}{\mathcal D} 
    \newcommand{\cP}{\mathcal P}
    \newcommand{\cC}{\mathcal C} 
    \newcommand{\cS}{\mathcal S}
    \newcommand{\EGC}{{\cal E}(\GC)}
    \newcommand{\cLGC}{\widetilde{L}_{an}\GC}
    \newcommand{\LGC}{{L}_{an}\GC} 
    \newcommand{\BQ}{\begin{eqnarray}}
    \newcommand{\EQ}{\end{eqnarray}} 
    \newcommand{\BQn}{\begin{eqnarray*}}
    \newcommand{\EQn}{\end{eqnarray*}} 
    \newcommand{\wtilde}{\widetilde}
    \newcommand{\Hol}{\mbox{Hol}} 
    \newcommand{\Hom}{\mbox{Hom}}
    \newcommand{\poly}{polynomial } 
    \newcommand{\polys}{polynomials }
    \newcommand{\pz}{\frac{\p}{\p z}} 
    \newcommand{\pzi}{\frac{\p}{\p z_i}}
    \newcommand{\edge}{\text{\raisebox{2.75pt}{\makebox[20pt][s]{\hrulefill}}}}
    \newcommand{\halfedge}{\text{\raisebox{2.75pt}{\makebox[10pt][s]{\hrulefill}}}}
    \newcommand{\n}{\notag}
    \newcommand{\C}{\mathbb C}
    \newcommand{\A}{\mathcal A}
    \newcommand{\Q}{\mathbb Q}
    \newcommand{\X}{X_1,\ldots,X_n}
\newcommand{\Xmi}{X,\hat i\,}
\newcommand{\Xmj}{X,\hat j\,}
\newcommand{\Xmij}{X,\hat i\,\hat j\,}
    \newcommand{\Z}{Z_1,\ldots,Z_n}
     \newcommand{\Zmi}{Z,\hat i}
    \newcommand{\pa}{\partial}
    \newcommand{\D}{D_1,\ldots,D_n}
    \newcommand{\Del}{\text{\raisebox{2.1pt}{$\bigtriangledown$}}}
    \newcommand{\La}{\triangle}
    \newcommand{\Hess}{\text{\rm Hess}}
    \newcommand{\mfs}{\mathfrak s}
    \newcommand{\mft}{\mathfrak t}
    \newcommand{\mfH}{\mathfrak H}
\newcommand{\F}{F_{1},\ldots,F_{n}}
\newcommand{\G}{G_{1},\ldots,G_{n}}
\newcommand{\GA}{\text{GA}}
\newcommand{\SA}{\text{SA}}
\newcommand{\MA}{\text{MA}}
\newcommand{\GL}{\text{GL}}
\newcommand{\SL}{\text{SL}}
\newcommand{\GE}{\text{GE}}
\newcommand{\PGL}{\text{PGL}}
\newcommand{\Tr}{\text{Tr}}
\newcommand{\Af}{\text{Af}}
\newcommand{\Bf}{\text{Bf}}
\newcommand{\W}{\text{W}}
\newcommand{\EA}{\text{EA}}
\newcommand{\BA}{\text{BA}}
\newcommand{\E}{\text{E}}
\newcommand{\B}{\text{B}}
\newcommand{\TA}{\text{TA}}
\newcommand{\Cr}{\text{Cr}}
\newcommand{\Di}{\text{D}}
\newcommand{\WT}{\text{WTA}}
\newcommand{\degr}{\text{deg}\,}
\newcommand{\supp}{\text{supp}\,}
\newcommand{\ideaala}{\mathfrak{a}}
\newcommand{\ideaalb}{\mathfrak{b}}
\newcommand{\ideaalc}{\mathfrak{c}}
\newcommand{\m}{\mathfrak{m}}
\newcommand{\id}{\text{id}}

\abstract{It is shown the the tame subgroup $\TA_3(\C)$ of the group $\GA_3(\C)$ of polynomials automorphisms of $\C^3$ can be realized as the product of three subgroups, amalgamated along pairwise intersections, in a manner that generalizes the well-known amalgamated free product structure of $\TA_2(\C)$ (which coincides with $\GA_2(\C)$ by Jung's Theorem).  The result follows from defining relations for $\TA_3(\C)$ given by U. U. Umirbaev.}

\section{Polynomial automorphism groups}

For a commutative ring $R$, we write $R^{[n]}$ for the polynomial ring $R[\X]$ in $n$ variables over $R$.  We will have occasion to refer to the subalgebra $R[X_1,\ldots,X_{i-1},X_{i+1},\ldots,X_n]$ for $i\in\{1,\ldots,n\}$, so we will use the shorter notation $R[\Xmi]$ to denote it.

The symbol $\GA_n(R)$ denotes the {\it general automorphism group, }by which we mean the automorphism group of $\text{Spec}\,R^{[n]}$ over $\text{Spec}\,R$.  As such, it is anti-isomorphic to the group of $R$-algebra automorphisms of $R^{[n]}$.   An element of $\GA_n(R)$ is represented by a vector $\varphi=(\F)\in(R^{[n]})^n$; we will consistently use Greek letters to denote automorphisms.  

The general linear group $\GL_n(R)$ is contained in $\GA_n(R)$ in an obvious way.  Another familiar subgroup is $\EA_n(R)$, the group generated by {\it elementary} automorphisms, i.e., those of the form $(X_1,\ldots,X_{i-1},X_i+f,X_{i+1},\ldots,X_n)$ for some $i\in\{1,\ldots,n\}$, $f\in R[\Xmi]$.

The subgroup of {\it tame automorphisms,} denoted $\TA_2(R)$, is the subgroup generated by $\GL_n(R)$ and $\EA_n(R)$.  Other subgroups of interest are the {\it affine }group $\Af_n(R)$, which is the group generated by $\GL_n(R)$ together with the {\it translations, }i.e., those automorphisms of the form $(X_1+a_1,\ldots,X_n+a_n)$ with $a_1,\ldots,a_n\in R$.

For $K$ a field, $\GA_n(K)$ is sometimes called the {\it affine Cremona group.} It sits naturally as a subgroup of the full Cremona group $\Cr_n(K)$, which is the group of birational automorphisms of affine (or projective) $n$-space.  The Jung-Van der Kulk Theorem (\cite{Jung},\cite{vdK}) states that $\TA_2(K)=\GA_2(K)$.  Shestakov and Umirbaev (\cite{S-U}) showed that $\TA_3(K)\ne\GA_3(K)$ when $K$ has characteristic zero.  This paper deals with the structure of $\TA_3(K)$, when char$(K)=0$, based of work of Umirbaev in \cite{UmirRel}.

\section{Amalgamated Products of Groups} \label{S1}

We begin with the definition of an amalgamation of groups.

\begin{defi}  Suppose we are given groups $A_i$ for each $i\in\{1,\cdots,n\}$  and for each $i,j\in \{1,\cdots,n\}$ with $i\ne j$ we have groups $B_{ij}=B_{ji}$ with injective homomorphisms $\varphi_{ij}:B_{ij}\to A_i$ which are compatible, meaning if $i,j,k$ are distinct then $\varphi_{ij}^{-1}(\varphi_{ik}(B_{ik}))=\varphi_{ji}^{-1}(\varphi_{jk}(B_{jk}))$ and on this group $\varphi_{ik}^{-1}\varphi_{ij} =\varphi_{jk}^{-1}\varphi_{ji}$.  This gives set-theoretic gluing data by which we can compatibly glue $A_i$ to $A_j$ along $B_{ij}$ via $\varphi_{ij}^{-1}\varphi_{ji}$ forming an amalgamated union $S$ of the sets $A_{1},\ldots,A_n$.  We then form the free group $\mathcal F$ on $S$, denoting the group operation on $\mathcal F$ by $*$.  For $i\in\{1,\cdots,n\}$ and $x,y\in A_i\subset S$, we let $r_{x,y}= x*y*(xy)^{-1}\in \mathcal F$ (where $xy$ is the product in $A_i$).  Finally we let $\mathcal G$ be the quotient of $\mathcal F$ by all the relations $r_{x,y}$.  The group $\mathcal G$ is called the {\it amalgamated product} of the groups $A_i$, $i\in\{1,\ldots,n\}$ along the groups $B_{ij}$, $i,j\in\{1,\ldots,n\}$.  There are natural group homomorphisms $\iota_i:A_i\to\mathcal G$ with $\iota_i\varphi_{ij}=\iota_j\varphi_{ji}$ on $B_{ij}$.
\end{defi}

The group $\mathcal G$ has the following universal property: Given a group $H$ and maps $\rho_i:A_i\to H$ for $i\in\{1,\ldots,n\}$ such that $\rho_i\varphi_{ij}=\rho_j\varphi_{ji}$ on $B_{ij}$ for all $i,j\in\{1,\ldots,n\}$, then there is a unique map $\Phi:\mathcal G\to H$ with $\rho_i\iota_i$ for all $i$.

When a group $\mathcal G$ is the amalgamation of two subgroups $A_1$ and $A_2$ along a common subgroup $B$, the two groups inject into the amalgamated product and a very strong factorization theorem holds.  Moreover the Bass-Serre tree theory of groups acting on trees (see \cite{SerreTrees}) provides a tree on which $\mathcal G$ acts without inversion, having a fundamental domain consisting of a single edge with its end vertices, the stabiliizers of the vertices being $A_1$ and $A_2$ and the stabilizer of the edge the common subgroup $B$.

Such theorems do not hold in general for amalgamations of three or more groups along pairwise intersections.  The groups $A_i$ may not map injectively into $\mathcal G$, and in fact $\mathcal G$ may be the trivial group when none of the groups $A_i$ are trivial, as the following example from \cite{Stall}  shows.

\begin{exam}
For $\{i,j,k\}=\{1,2,3\}$ let $B_{ij}$ be the infinite cyclic group generated by $b_k$.  Let
\begin{align*}
A_1&=\langle\,b_2,b_3\,|\,b_2b_3b_2^{-1}=b_3^2\,\rangle\\
A_2&=\langle\,b_3,b_1\,|\,b_3b_1b_3^{-1}=b_1^2\,\rangle\\
A_3&=\langle\,b_1,b_2\,|\,b_1b_2b_1^{-1}=b_2^2\,\rangle
\end{align*}
Then $B_{ij}$ is a common subgroup of $A_i$ and $A_j$ and we can form the amalgamation $\mathcal G$ of the groups $A_i$ along the groups $B_{ij}$.  It can be shown that in this case $\mathcal G$ is the trivial group.
\end{exam}

Whether such amalgamation data gives rise to the group acting on a simplicial complex is not easy to detect (see, for example, \cite{Stall}, \cite{Haef}, and \cite{Cor}).  It occurs precisely when each of the groups $A_{ij}$ maps injectively to $\mathcal G$, and in this situation, the amalgamated union $S$ maps injectively to $\mathcal G$ as well.  The $n$-simplex of groups arising from this data is called {\it developable} by Haefliger  (\cite{Haef}) in case of this occurance.

However, if the groups $A_i$ are subgroups of a given group $G$ and if we take $B_{ij}$ to be $A_i\cap A_j$ and $\varphi_{ij}$ the inclusion map within $G$, then clearly there exists a homomorphism $\Phi:\mathcal G\to G$ restricting to the identity on each $A_i$, which shows that in this case the amalgamated union $S$ maps injectively to $\mathcal G$.  The map $\Phi$ will be surjective precisely when $G$ is generated by the subgroups $A_1,\ldots, A_n$.  If $\Phi$ is an isomorphism, then the structure of $\mathcal G$ arises from the action of $G$ on an $n$-dimensional simply connected simplicial complex, with a single simplex serving as a fundamental domain.

Automorphism groups of various kinds can be realized as amalgamations of groups.  Some examples are given below.  

\begin{exam}
$\SL_2(\mathbb Z)=(\mathbb Z/4\mathbb Z)*_{\mathbb Z/2\mathbb Z}(\mathbb Z/6\mathbb Z)$ acts on the upper half plane.  The generator of $\mathbb Z/4\mathbb Z)$ and $\mathbb Z/6\mathbb Z$ can be taken to be $\left(\begin{smallmatrix}0&1\\-1&0\end{smallmatrix}\right)$ and $\left( \begin{smallmatrix}1&-1\\1&0\end{smallmatrix}\right)$, respectively.  Here the translates of the circular arc $z=e^{i\theta}$ with $\pi/3\le\theta\le\pi/2$ form a tree with this arc as a fundamental domain, and this is the tree given by the Bass-Serre theory.
\end{exam}

\begin{exam}[Nagao's Theorem \cite{Nagao}]
For $K$ a field we have $\GL_2(K[X])=\GL_2(K)*_{\B_2(K)}\B_2(K[X])$, with $\B_2$ denoting the lower triangular group.  This structure can be realized via the Bass-Serre theory by the action of $\GL_2(K[X])$  on a tree whose vertices are $\mathcal O$-lattices in the rank two vector space over $K(X)$. where $\mathcal O$ is the DVR of $K(X)$ with uniformizing parameter $1/X$ (see \cite{SerreTrees}).  Here the fundamental domain is not just a single edge, but an edge connected to a ``directed geodesic.''  
\end{exam}

\begin{exam}[Jung-Van der Kulk Theorem \cite{Jung},\cite{vdK}]\label{JVdK}
For $K$ a field, the group $\GA_2(K)$ of polynomial automorphisms of the affine plane has he structure $\GA_2(K))=\Af_2(K)*_{\Bf_2(K)}\BA_2(K)$.  Here $\BA_2$ is the group of automorphisms of the form $(X_1+\alpha, X_2+f(X_1))$,  $f(X_1)\in K[X_1]$ and $\Bf_2(K)=\Af_2(K)\cap\BA_2(K)$.  Again, this structure arises from the action of $\GA_2(K)$ on a tree whose vertices are certain complete algebraic surfaces realized as collections of local rings (``models") inside the function field $K(X_1,X_2)$ (see \cite{W8}).
\end{exam}

\begin{exam}
The full Cremona group $\Cr_2(K)$ over an algebraically closed field $K$ is the amalgamation of three groups: the automorphism group of $\mathbb P_K^2$ (which is $\PGL_2(K)$), the automorphism group of $\mathbb P_K^1\times\mathbb P_K^1$, and thirdly the $K$-automorphism group of $\mathbb P_L^1$ where $L=K(t)$, with $t$ transcendental over $K$.  There is a naturally realizable simplicial complex of triangles $\mathcal C$ on which  $\Cr_2(K)$ acts which yields this structure and also contains the tree of Example \ref{JVdK} with the action of $\GA_2(K)$ being the restriction of the action of $\Cr_2(K)$ on $\mathcal C$.  See \cite{W8} for details.
\end{exam}

\section{Polynomial Automorphisms in Dimension Three} \label{dim3}

A major breakthrough came in 2004 when Shestakov and Umirbaev showed that the automorphism group $\GA_3(K)$ properly contains the tame subgroup $\TA_3(K)$ when $K$ is a field of characteristic zero (\cite{S-U}).  Specifically they showed that the automorphism
$$\left(X+Z(YZ+X^2),Y-2X(YZ+X^2)-X(YZ+X^2)^2,Z\right)$$
is not tame, resolving a conjecture of Nagata from 1972.  The group $\GA_3(K)$ remains a mystery, as no describable set of generators has been given.

However, the tame subgroup $\TA_3(K)$ is (by definition) generated by the elementary and linear automorphisms, which are familiar.  Moreover a set of generating relations has been given by Umirbaev in \cite{UmirRel}.  This paper will show that an amalgamated product structure for $\TA_3(K)$ results from Umirbaev's relations.   We begin by presenting those results.

For $\varphi=(\F)\in\GA_n(K)$ and $f\in K^{[n]}$ we write $f(\varphi)$ for $f(\F)$.  This defines an action (on the right) of $\GA_n(K)$ on $K^{[n]}$.

For $i\in\{1,\ldots,n\}$, $\alpha\in K$, and $f\in K[\Xmi]$, consider the automorphism
\begin{equation}\label{sigma}
\sigma_{i,\alpha,f}=(X_1,\ldots,X_{n-1},\alpha X_i+f,X_{i+1},\ldots,X_n)\,,
\end{equation}
which is easily seen to be tame.  Given $k,\ell\in\{1,\ldots,n\}$, $k\ne\ell$, we define a tame automorphism $\tau_{k.\ell}$ by
\begin{equation}\label{tau}
\tau_{k.\ell}=\sigma_{\ell,1,X_k}\sigma_{k,1,-X_\ell}\sigma_{\ell,-1,X_k}
\end{equation}
A simple calculation show that $\tau_{k,\ell}$ is the transposition switching the $X_k$ and $X_\ell$ coordinates. 

One can check directly that
\begin{equation}\label{rel1}
\sigma_{i,\alpha,f}\sigma_{i,\beta,g}=\sigma_{i,\alpha\beta,f+\alpha g}\,.
\end{equation}
Also, if $i,j\in\{1,\ldots,n\}$, $i\ne j$, and if $f\in  K[\Xmi]\cap  K[\Xmj]$, $g\in K[\Xmj]$, then
\begin{equation}\label{rel2}
{\sigma_{i,\alpha,f}}^{-1}\sigma_{j,\beta,g}\sigma_{i,\alpha,f}=\sigma_{j\beta,g(\sigma_{i,\alpha,f})}
\end{equation}
It follows that if $g\in  K[\Xmi]\cap  K[\Xmj]$ then $\sigma_{i,\alpha,f}$ and $\sigma_{i,\beta,g}$ commute.

Let $k,\ell\in\{1,\ldots,n\}$, $k\ne\ell$.  For $i\in\{1,\ldots,n\}$ let $j$ be the image of $i$ under the permutation which switches $k$ and $\ell$, in other words, the element of $\{1,\ldots,n\}$ for which $X_j=X_i(\tau_{k,\ell})$.  Then we have
\begin{equation}\label{rel3}
\tau_{k,\ell}\sigma_{i,\alpha,f}\tau_{k,\ell}=\sigma_{j,\alpha,f(\tau_{k,\ell})}
\end{equation}

Theorem 4.1 of \cite{UmirRel} asserts the following.
\begin{theo}[Umirbaev]\label{umir} Let $K$ be a field of characteristic zero. The relations (\ref{rel1}), (\ref{rel2}), and (\ref{rel3}) are defining relations for $\TA_3(K)$ with respect to the generators $\sigma_{i,\alpha,f}$ defined in (\ref{sigma}).  Here $\tau_{k,\ell}$ in (\ref{rel3}) is defined formally in terms of these generators by (\ref{tau}).
\end{theo}

\noindent This will be the key tool in the proof of Theorem \ref{main}, which is the main result of this paper.

\section{Subgroups of Interest} \label{S4}

For $i\in\{1,\ldots,n\}$, let $V_i$ be the sub-vector space of $K^{[n]}$ generated by $K$ and the variables $X_1,\ldots,X_i$, i.e.,
\begin{equation}\label{vecsp}
V_i=K\oplus KX_1\oplus\cdots\oplus KX_i\,.
\end{equation}
Then $H_i$ is defined to be the stabilizer of $V_i$ in $\GA_n(K)$ via the action defined in \S\ref{dim3}, i.e., 
\begin{equation}\label{Hdef}
H_i=\{\varphi\in\GA_n(K)\,|\, f(V_i)=V_i\}\,.
\end{equation}  
Note that $H_n$ is the affine group $\Af_n(K)$.  More generally, the subgroup of $H_i$ that fixes each of the variables $X_{i+1},\ldots,X_n$ can be identified with $\Af_i(K)$.   In fact, $H_i$ retracts onto $\Af_i(K)$ via the map $\varphi=(\F)\mapsto(F_1,\ldots,F_i)$, and the kernel of this retraction is the subgroup of $H_i$ consisting of the elements that fix each of the variables $X_1,\ldots,X_i$, which is $\GA_{n-1}(K[X_1,\ldots,X_i])$.  Thus $H_i$ has the semidirect product structure
\begin{equation}\label{Hi}
H_i=\Af_i(K)\ltimes\GA_{n-i}(K[X_1,\ldots,X_i])
\end{equation}
(where, for $i=n$, we read this as $H_n=\Af_n(K)$).  These subgroups are defined in \cite{F-R}, p.\ 23, where it is conjectured that together they generate $\GA_n(K)$ (Conjecture 14.1) and that (whether or not that conjecture is true) the subgroup generated by $H_1,\ldots,H_n$ is the amalgamated product of these groups along pairwise intersections (Conjecture 14.2).  It should be noted that Freudenburg produced an example (see \cite{Fr1}, p.\ 121) of an automorphism in $GA_3(K)$ which has not been shown to lie in this subgroup.\footnote{This example is also of interest because it has not been shown to be stably tame.}

Furthermore the groups $\widetilde H_i$ are defined by
$$\widetilde H_i=H_i\cap\TA_n(K)\,,$$
which are easily seen to generate $\TA_n(K)$.  We can surmise from (\ref{Hi}) that 
\begin{equation}\label{widetildeHi}\widetilde H_i\supseteq\Af_i(K)\ltimes\TA_{n-i}(K[X_1,\ldots,X_i])\,.
\end{equation}  For $i=n$ equality holds trivially and we have $\widetilde H_n=H_n$, both being equal to $\Af_n(K)$.  For $i=n-1$ it is also easily seen that equality holds in (\ref{widetildeHi})  and moreover we have $\widetilde H_{n-1}=H_{n-1}$ since $\TA_{1}$ and $\GA_{1}$ coincide over an integral domain (even a reduced ring).  

There is one other case where the containment of (\ref{widetildeHi}) is known to be an equality.  Namely, for $n=3$ and $K$ of characteristic zero we have $\widetilde H_1=Af_1(K)\ltimes\TA_2(K[X_1])$.  This follows from Corollary 10 of \cite{S-U}, a very deep result asserting that in $\GA_3(K)$ we have $$\GA_2(K[X_1])\cap\TA_3(K)=\TA_2(K[X_1])\,.$$  This together with the known proper containment $\TA_2(K[X_1])\subsetneqq\GA_2(K[X_1])$ tells us that $\widetilde H_1\subsetneqq H_1$ for $n=3$.  It is not known whether $\widetilde H_1\subsetneqq H_1$ when $n>3$.

It is conjectured that $\TA_n(K)$ is the amalgamated product of the subgroups  $\widetilde H_1,\ldots,\widetilde H_n$ along pairwise intersections (\cite{F-R}, Conjecture 14.3).  The main result of this paper is that this conjecture is true for $n=3$ and $K$ a field of characteristic zero.  In light of the above observations, for $n=3$ we have $\widetilde H_2=H_2$ and $\widetilde H_3=H_3$ (but \underbar{not} $\widetilde H_1=H_1$), so this can be stated as:

\begin{theo}\label{main} For $K$ a field of characteristic zero, $\TA_3(K)$ is the amalgamated product of the three groups $\widetilde H_1,H_2,H_3$ along their pairwise intersections.
\end{theo}

\noindent This will be proved in the next section.  %We first make some observations about these subgroups.  As noted above,
  
%This means we have a triangle of groups, in Stallings' sense.

\section{Proof of Theorem \ref{main}} 

The main tool in the proof is Theorem \ref{umir}, which asserts that $\TA_2(K)$ is generated by the elements $\sigma_{i,\alpha,f}$ as defined in (\ref{sigma}) subject to the relations (\ref{rel1}), (\ref{rel2}), and (\ref{rel3}).  

Let $\mathcal F$ be the free group generated by the formal symbols $[\sigma_{i,\alpha,f}]$, with $i\in\{1,\ldots,n\}$, $f\in K[\Xmi]$.  Accordingly, we rewrite the relations (\ref{rel1}), (\ref{rel2}), and (\ref{rel3}) replacing each $\sigma$ by its corresponding formal symbol $[\sigma]$:
\begin{gather}
[\sigma_{i,\alpha,f}][\sigma_{i,\beta,g}]=[\sigma_{i,\alpha\beta,f+\alpha g}]\label{rel1a}\tag{R1}\\
{[\sigma_{i,\alpha,f}]}^{-1}[\sigma_{j,\beta,g}][\sigma_{i,\alpha,f}]=[\sigma_{j,\beta,g(\sigma_{i,\alpha,f})}]\label{rel2a}\tag{R2}\\
[\tau_{k,\ell}][\sigma_{i,\alpha,f}][\tau_{k,\ell}]=[\sigma_{j,\alpha,f(\tau_{k,\ell})}]\,,\label{rel3a}\tag{R3}
\end{gather}
where, in (\ref{rel2a}), $i\ne j$,  $f\in  K[\Xmi]\cap  K[\Xmj]$, $g\in K[\Xmj]$, and, in (\ref{rel3a}), $k\ne\ell$,  $j$ is the image of $i$ under the permutation which switches $k$ and $\ell$, and
\begin{equation}\label{[tau]}
[\tau_{k,\ell}]=[\sigma_{\ell,1,X_k}][\sigma_{k,1,-X_\ell}][\sigma_{\ell,-1,X_k}]
\end{equation}
(after (\ref{tau})).  Let $\mathcal N$ be the normal subgroup of $\mathcal F$ generated by (\ref{rel1a}), (\ref{rel2a}), and (\ref{rel3a}).

Theorem \ref{umir} says that the homomorphism from $\mathcal F$ to $\TA_3(K)$ sending $[\sigma_{i,\alpha,f}]$ to $\sigma_{i,\alpha,f}$ induces an isomorphism  
\begin{equation}\label{isom}
\mathcal F/\mathcal N\overset{\cong}{\longrightarrow}\TA_3(K)\,.
\end{equation}

Let $\mathfrak G$ be the amalgamated product of $\widetilde H_1, H_2, H_3$ along their pairwise intersections.  The inclusions of $\widetilde H_1, H_2, H_3$ in $\TA_3(K)$ induce a group homomorphism $\Phi:\mathfrak G\to\TA_3(K)$ which is surjective since the three subgroups generate $\TA_3(K)$ (in fact any two of them generate).  We will define a group homomorphism from $\TA_3(K)$ to $\mathfrak G$ using the isomorphism (\ref{isom}) and show that it is inverse to $\Phi$, thus proving the theorem.  

We first define a homomorphism $\widehat\Psi:\mathcal F\to\mathfrak G$, which is accomplished by specifying the images of the free generators $[\sigma_{i,\alpha,f}]$.  According to the discussion in Section \ref{S1}, $\mathfrak G$ contains the amalgamated union of $\widetilde H_1, H_2, H_3$ as does $TA_3(K)$, with $\Phi$ restricting to the identity map on this set.  Let us denote by $\widetilde{\mfH}_1,{\mfH}_2, {\mfH}_3$ the isomorphic copies of $\widetilde H_1, H_2, H_3$, respectively, that lie inside $\mathfrak G$.  It is important to keep in mind that $\widetilde{\mfH}_1\cup{\mfH}_2\cup{\mfH}_3$ maps bijectively to $\widetilde H_1\cup H_2\cup H_3$ via $\Phi$.

Note that if $i=2$ or $i=3$ then $\sigma_{i,\alpha,f}$ lies in $\widetilde H_1$ and if $\text{deg}\,f\le1$ then $\sigma_{i,\alpha,f}$ lies in $H_3$, so in each of these cases $\sigma_{i,\alpha,f}$ can be viewed as an element of the union $\widetilde{\mfH}_1\cup{\mfH}_2\cup{\mfH}_3\subset\mathfrak G$.  To avoid confusion, we will denote these elements of $\mathfrak G$ by $\mfs_{i,\alpha,f}$.  Thus it makes sense to make the assignments
\begin{equation}\label{assign23}
\widehat\Psi([\sigma_{i,\alpha,f}])=\mfs_{i,\alpha,f}\in\widetilde{\mfH}_1\quad\text{for }i=2,3\,.
\end{equation}
\begin{equation}\label{assign1}
\widehat\Psi([\sigma_{1,\alpha,f}])=\mfs_{1,\alpha,f}\in \mfH_3\quad\text{for deg}\,f\le1\,.
\end{equation}
Since the factors of (\ref{[tau]}) involve only polynomials of degree $\le1$,  $\widehat\Psi([\tau_{k,\ell}])$ is defined by aplying $\widehat\Psi$ to those factors using (\ref{assign23}) and(\ref {assign1}) above.  We will denote the resulting element of $\mathfrak G$ by $\mft_{k,\ell}\,$.  Thus:
\begin{equation}\label{Phi[tau]}
\widehat\Psi([\tau_{k,\ell}])=\mft_{k,\ell}=\mfs_{\ell,1,X_k}\,\mfs_{k,1,-X_\ell}\,\mfs_{\ell,-1,X_k}\,.
\end{equation}
and this is the just the permutation in $\mfH_3\cong\Af_3(K)$ that switches $k$ and $\ell$.  It remains to define $\widehat\Psi([\sigma_{1,\alpha,f}])$ for arbitrary $f\in K[X_2,X_3]$.   This we do as follows:
\begin{equation}\label{assignf}
\widehat\Psi\left([\sigma_{1,\alpha,f}]\right)=\mft_{1,3}\,\mfs_{3,\alpha,f(X_2,X_1)}\,\mft_{1,3}\,.
\end{equation}
The reader will easily verify that this assignment coincides with (\ref{assign1}) in the case $\text{deg}\,f\le1$, since both occur in $\mfH_3$.

Thus we have defined $\widehat\Psi:\mathcal F\to\mathfrak G$, and we must now show that the subgroup $\mathcal N$ lies in the kernel of $\widehat\Psi$, i.e, that equations (\ref{rel1a}), (\ref{rel2a}), and (\ref{rel3a}) hold replacing $\sigma$ by $\mfs$ and $\tau$ by $\mft$.  This gets a bit tedious because of the asymmetry in the definitions of $\widehat\Psi\left([\sigma_{i,\alpha,f}]\right)$ depending on $i$.

We begin with (\ref{rel1a}).  Note that if $i=2$ or $i=3$, then according to (\ref{assign23}), this amounts to showing that 
\begin{equation}\label{i23}
\mfs_{i,\alpha,f}\,\mfs_{i,\beta,g}=\mfs_{i,\alpha\beta,f+\alpha g}\quad\text{for $i=2,3$.}
\end{equation}  
But this is a relation that takes place in $\widetilde{\mfH}_1$, so it holds in $\mathfrak G$.  For $i=1$ we must use (\ref{assignf}).  For $f,g\in K[X_2,X_3]$ we have
\begin{align*}
\widehat\Psi([\sigma_{1,\alpha,f}])\widehat\Psi([\sigma_{1,\beta,g}])&=(\mft_{1,3}\,\mfs_{3,\alpha,f(X_2,X_1)}\,\mft_{1,3})(\mft_{1,3}\,\mfs_{3,\beta,g(X_2,X_1)}\,\mft_{1,3})\\
&=\mft_{1,3}\,\mfs_{3,\alpha,f(X_2,X_1)}\,\mfs_{3,\beta,g(X_2,X_1)}\,\mft_{1,3}\\
&\qquad\qquad\qquad\qquad\quad\text{(since $\mft_{1,3}^2=1$ in $\mfH_3$)}\\
&=\mft_{1,3}\,\mfs_{3,\alpha\beta,f+\alpha g}\,\mft_{1,3}\qquad\text{by (\ref{i23})}\\
&=\widehat\Psi([\sigma_{1,\alpha\beta,f+\alpha g}])\,,\qquad\text{by (\ref{assignf})}
\end{align*}
completing the proof that the relation (\ref{rel1a}) is respected by $\widehat\Psi$.

We now address (\ref{rel2a}).  If $\{i,j\}=\{2,3\}$ we must show, again appealing to (\ref{assign23}), that
\begin{equation}\label{rel2a23}
\mfs_{i,\alpha,f}^{-1}\,\mfs_{j,\beta,g}\,\mfs_{i,\alpha,f}=\mfs_{j\beta,g(\sigma_{i,\alpha,f})}\quad\text{for $i=2,3$.}
\end{equation}
But, again, this is a relation that holds in $\widetilde{\mfH}_1$, hence in $\mathfrak G$.  

We now consider the case $i=1$, $j=3$.  We will use the following basic permutation relation, which holds in the symmetric group  $\mathfrak S_3\subset\mfH_3$ (hence it holds in $\mathfrak G$) for $\{k,\ell,m\}=\{1,2,3\}$:
\begin{equation}\label{permrel}
\mft_{k,\ell}=\mft_{k,m}\,\mft_{m,\ell}\,\mft_{k,m}\,.
\end{equation}

In the equations below the underbrace indicates what will be replaced in the next line; the overbrace in the next line marks the equivalent expression that has been substituted.  

For $f\in K[X_2]$ and $g\in K[X_1,X_2]$,
\begin{align*}
\widehat\Psi&\left([\sigma_{1,\alpha,f(X_2)}]^{-1}\,[\sigma_{3,\beta,g(X_1,X_2)}]\,[\sigma_{1,\alpha,f(X_2)}]\right)\\
&=\widehat\Psi\left({[\sigma_{1,\alpha,f(X_2)}]}\right)^{-1}\,\widehat\Psi\left([\sigma_{3,\beta,g(X_1,X_2)}]\right)\,\widehat\Psi\left([\sigma_{1,\alpha,f(X_2)}]\right)\\
&=\left(\underbrace{\mft_{1,3}}\,{\mfs_{3,\alpha,f(X_2)}}^{-1}\,\underbrace{\mft_{1,3}}\right)\,\left(\mfs_{3,\beta,g(X_1,X_2)}\right)\,\left(\underbrace{\mft_{1,3}}\,\mfs_{3,\alpha,f(X_2)}\,\underbrace{\mft_{1,3}}\right)\qquad\text{by (\ref{assign23}) and (\ref{assignf})}\\
\intertext{Applying (\ref{permrel}) to $\mft_{1,3}\,$:}
&=\lefteqn{\overbrace{\phantom{\mft_{1,2}\,\mft_{2,3}\,\mft_{1,2}}}}\,\mft_{1,2}\,\mft_{2,3}\,
\lefteqn{\underbrace{\phantom{\mft_{1,2}\,{\mfs_{3,\alpha,f(X_2)}}^{-1}\,\mft_{1,2}}}}
{\mft_{1,2}\,{\mfs_{3,\alpha,f(X_2)}}^{-1}\,
\lefteqn{\overbrace{\phantom{\mft_{1,2}\,\mft_{2,3}\,\mft_{1,2}}}}
\mft_{1,2}}\,\mft_{2,3}\,\mft_{1,2}\,\mfs_{3,\beta,g(X_1,X_2)}\,
\lefteqn{\overbrace{\phantom{\mft_{1,2}\,\mft_{2,3}\,\mft_{1,2}}}}
\mft_{1,2}\,\mft_{2,3}\,
\lefteqn{\underbrace{\phantom{\mft_{1,2}\,{\mfs_{3,\alpha,f(X_2)}}\,\mft_{1,2}}}}
\mft_{1,2}\,{\mfs_{3,\alpha,f(X_2)}}\,\overbrace{\mft_{1,2}\,\mft_{2,3}\,\mft_{1,2}}\\
\intertext{Using the relation $\mft_{1,2}\,{\mfs_{3,\alpha,f(X_2)}}\,\mft_{1,2}=\mfs_{3,\alpha,f(X_1)}$ from $\mfH_2\,$:}
&=\mft_{1,2}\,\mft_{2,3}\,\overbrace{{\mfs_{3,\alpha,f(X_1)}}^{-1}}\,\mft_{2,3}\,\underbrace{\mft_{1,2}\,\mfs_{3,\beta,g(X_1,X_2)}\,
\mft_{1,2}}\,\mft_{2,3}\,\overbrace{{\mfs_{3,\alpha,f(X_1)}}}\,\mft_{2,3}\,\mft_{1,2}
\intertext{Using the relation $\mft_{1,2}\,{\mfs_{3,\beta,g(X_1,X_2)}}\,\mft_{1,2}=\mfs_{3,\beta,g(X_2,X_1)}$ from $\mfH_2\,$:}
&=\mft_{1,2}\,\mft_{2,3}\,{\mfs_{3,\alpha,f(X_1)}}^{-1}\,\underbrace{\mft_{2,3}\,\overbrace{\mfs_{3,\beta,g(X_2,X_1)}}\,
\mft_{2,3}}\,{\mfs_{3,\alpha,f(X_1)}}\,\mft_{2,3}\,\mft_{1,2}\\
\intertext{Using the relation $\mft_{2,3}\,{\mfs_{3,\beta,g(X_2,X_1)}}\,\mft_{2,3}=\mfs_{2,\beta,g(X_3,X_1)}$ from $\widetilde{\mfH}_1\,$:}
&=\mft_{1,2}\,\mft_{2,3}\,\underbrace{{\mfs_{3,\alpha,f(X_1)}}^{-1}\,\overbrace{\mfs_{2,\beta,g(X_3,X_1)}}\,
{\mfs_{3,\alpha,f(X_1)}}}\,\mft_{2,3}\,\mft_{1,2}\\
\intertext{Applying (\ref{rel2a23}):}
&=\mft_{1,2}\,\underbrace{\mft_{2,3}\,\overbrace{\mfs_{2,\beta,g(\alpha X_3+f(X_1),X_1)}}\,\mft_{2,3}}\,\mft_{1,2}\\
\intertext{Using the relation $\mft_{2,3}\,\mfs_{2,\beta,g(\alpha X_3+f(X_1),X_1)}\,\mft_{2,3}=\mfs_{3,\beta,g(\alpha X_2+f(X_1),X_1)}$ from $\widetilde{\mfH}_1\,$:}
&=\underbrace{\mft_{1,2}\,\overbrace{\mfs_{3,\beta,g(\alpha X_2+f(X_1),X_1)}}\,\mft_{1,2}}\\
&=\mfs_{3,\beta,g(\alpha X_1+f(X_2),X_2)}\qquad\text{from $\mfH_2$}\\
&=\widehat\Psi\left([\sigma_{3,\beta,g(\alpha X_1+f(X_2),X_2)}]\right)\\
&=\widehat\Psi\left([\sigma_{3,\beta,g(\sigma_{1,\alpha,f(X_2)})}]\right)\,,
\end{align*}
which accomplishes our goal.  

%\[
%0+\lefteqn{\overbrace{\phantom{1+2+3}}}1+
%\underbrace{2+3+\overbrace{x+y+z}+
%\lefteqn{\overbrace{\phantom{4+5}}}4}+5
%\]

Now let $i=3$, $j=1$.  For $f\in K[X_2]$ and $g\in K[X_2,X_3]$,
\allowdisplaybreaks{\begin{align*}
\widehat\Psi&\left([\sigma_{3,\alpha,f(X_2)}]^{-1}\,[\sigma_{1,\beta,g(X_2,X_3)}]\,[\sigma_{3,\alpha,f(X_2)}]\right)\\
&=\widehat\Psi\left({[\sigma_{3,\alpha,f(X_2)}]}\right)^{-1}\,\widehat\Psi\left([\sigma_{1,\beta,g(X_2,X_3)}]\right)\,\widehat\Psi\left([\sigma_{3,\alpha,f(X_2)}]\right)\\
&={\mfs_{3,\alpha,f(X_2)}}^{-1}\,\mft_{1,3}\,\mfs_{3,\beta,g(X_2,X_1)}\,\mft_{1,3}\,\mfs_{3,\alpha,f(X_2)}\qquad\text{by (\ref{assign23}) and (\ref{assignf})}\\
&=\mft_{1,3}\,\underbrace{\mft_{1,3}}\,{\mfs_{3,\alpha,f(X_2)}}^{-1}\,\underbrace{\mft_{1,3}}\,\mfs_{3,\beta,g(X_2,X_1)}\,\underbrace{\mft_{1,3}}\,\mfs_{3,\alpha,f(X_2)}\,\underbrace{\mft_{1,3}}\,\mft_{1,3}\qquad\text{since ${\mft_{1,3}}^2=1$}\\
\intertext{Applying (\ref{permrel}):}
&=\mft_{1,3}\,\lefteqn{\overbrace{\phantom{\mft_{1,2}\,\mft_{2,3}\,\mft_{1,2}}}}
\mft_{1,2}\,\mft_{2,3}\,\lefteqn{\underbrace{\phantom{\mft_{1,2}\,{\mfs_{3,\alpha,f(X_2)}}^{-1}\,\mft_{1,2}}}}
\mft_{1,2}\,{\mfs_{3,\alpha,f(X_2)}}^{-1}\,
\lefteqn{\overbrace{\phantom{\mft_{1,2}\,\mft_{2,3}\,\mft_{1,2}}}}
\mft_{1,2}\,\mft_{2,3}\,\mft_{1,2}\,\mfs_{3,\beta,g(X_2,X_1)}\\
&\hskip2in\lefteqn{\overbrace{\phantom{\mft_{1,2}\,\mft_{2,3}\,\mft_{1,2}}}}
\mft_{1,2}\,\mft_{2,3}\,\lefteqn{\underbrace{\phantom{\mft_{1,2}\,\mfs_{3,\alpha,f(X_2)}\,\mft_{1,2}}}}
\mft_{1,2}\,\mfs_{3,\alpha,f(X_2)}\,\overbrace{\mft_{1,2}\,\mft_{2,3}\,\mft_{1,2}}\mft_{1,3}\\
\intertext{Using the relation $\mft_{1,2}\,\mfs_{3,\alpha,f(X_2)}\,\mft_{1,2}=\mfs_{3,\alpha,f(X_1)}$ in $\mfH_2$\,:}
&=\mft_{1,3}\,\mft_{1,2}\,\underbrace{\mft_{2,3}\,\overbrace{{\mfs_{3,\alpha,f(X_1)}}^{-1}}\,\mft_{2,3}}\,\mft_{1,2}\,\mfs_{3,\beta,g(X_2,X_1)}\,\mft_{1,2}\,\underbrace{\mft_{2,3}\,\overbrace{\mfs_{3,\alpha,f(X_1)}}\,\mft_{2,3}}\,\mft_{1,2}\,\mft_{1,3}\\
\intertext{Using the relation $\mft_{2,3}\,\mfs_{3,\alpha,f(X_1)}\,\mft_{2,3}=\mfs_{2,\alpha,f(X_1)}$ in $\widetilde{\mfH}_1\,$:}
&=\mft_{1,3}\,\mft_{1,2}\,\overbrace{{\mfs_{2,\alpha,f(X_1)}}^{-1}}\,\underbrace{\mft_{1,2}\,\mfs_{3,\beta,g(X_2,X_1)}\,\mft_{1,2}}\,\overbrace{\mfs_{2,\alpha,f(X_1)}}\,\mft_{1,2}\,\mft_{1,3}\\
\intertext{Using the relation $\mft_{1,2}\,\mfs_{3,\beta,g(X_2,X_1)}\,\mft_{1,2}=\mfs_{3,\beta,g(X_1,X_2)}$ in $\mfH_2\,$:}
&=\mft_{1,3}\,\mft_{1,2}\,\underbrace{{\mfs_{2,\alpha,f(X_1)}}^{-1}\,\overbrace{\mfs_{3,\beta,g(X_1,X_2)}}\,\mfs_{2,\alpha,f(X_1)}}\,\mft_{1,2}\,\mft_{1,3}\\
\intertext{Applying (\ref{rel2a23}):}
&=\mft_{1,3}\,\underbrace{\mft_{1,2}\,\overbrace{\mfs_{3,\beta,g(X_1,\alpha X_2+f(X_1))}}\,\mft_{1,2}}\,\mft_{1,3}\\
\intertext{Using the relation $\mft_{1,2}\,\mfs_{3,\beta,g(X_1,\alpha X_2+f(X_1))}\,\mft_{1,2}=\mfs_{3,\beta,g(X_2,\alpha X_1+f(X_2))}$ in $\mfH_2\,$:}
&=\mft_{1,3}\,\overbrace{\mfs_{3,\beta,g(X_2,\alpha X_1+f(X_2))}}\,\mft_{1,3}\\
&=\widehat\Psi\left([\sigma_{1,\beta,g(X_2,\alpha X_1+f(X_2))}]\right)\\
	&=\widehat\Psi\left([\sigma_{1,\beta,g(\sigma_{3,\alpha,f(X_2)})}]\right)\quad\text{by (\ref{assignf}),}
\end{align*}}
as desired.

The two cases $\{i,j\}=\{1,3\}$ will employ the equality
\begin{equation}\label{special12}
\mft_{1,3}\,\mfs_{2,\beta,g(X_1,X_3)}\,\mft_{1,3}=\mfs_{2,\beta,g(X_3,X_1)}\,,
\end{equation} 
which arises by conjugating the $\mfH_2$ identity $\mft_{1,2}\,\mfs_{3,\beta,g(X_1,X_2)}\,\mft_{1,2}=\mfs_{3,\beta,g(X_2,X_1)}$ by $\mft_{2,3}$, evoking the $\widetilde{\mfH}_1$ identity $\mft_{2,3}\,\mfs_{3,\beta,g(X_1,X_2)}\,\mft_{2,3}=\mfs_{2,\beta,g(X_1,X_3)}$ and the $\mfH_3$ identity $\mft_{2,3}\,\mft_{1,2}\,\mft_{2,3}=\mft_{1,3}$.

For $i=1, j=2$ we have
\begin{align*}
\widehat\Psi&\left([\sigma_{1,\alpha,f(X_3)}]^{-1}\,[\sigma_{2,\beta,g(X_1,X_3)}]\,[\sigma_{1,\alpha,f(X_3)}]\right)\\
&=\widehat\Psi\left({[\sigma_{1,\alpha,f(X_3)}]}\right)^{-1}\,\widehat\Psi\left([\sigma_{2,\beta,g(X_1,X_3)}]\right)\,\widehat\Psi\left([\sigma_{1,\alpha,f(X_3)}]\right)\\
&=\mft_{1,3}\,{\mfs_{3,\alpha,f(X_1)}}^{-1}\,\underbrace{\mft_{1,3}\,\mfs_{2,\beta,g(X_1,X_3)}\,\mft_{1,3}}\,\mfs_{3,\alpha,f(X_1)}\,\mft_{1,3}\qquad\text{by (\ref{assignf})}\\
&=\mft_{1,3}\,\underbrace{{\mfs_{3,\alpha,f(X_1)}}^{-1}\,\overbrace{\mfs_{2,\beta,g(X_3,X_1)}}\,\mfs_{3,\alpha,f(X_1)}}\,\mft_{1,3}\qquad\text{by (\ref{special12})}\\
&=\underbrace{\mft_{1,3}\,\overbrace{\mfs_{2,\beta,g(\alpha X_3+f(X_1),X_1)}}\,\mft_{1,3}}\qquad\text{by (\ref{rel2a23})}\\
&=\mfs_{2,\beta,g(\alpha X_1+f(X_3),X_3)}\qquad\text{by (\ref{special12})}\\
&=\mfs_{2,\beta,g(\sigma_{1,\alpha,f(X_3)})}\\
&=\widehat\Psi\left([\sigma_{2,\beta,g(\sigma_{1,\alpha,f(X_3)})}]\right)\,.
\end{align*}

The case $i=2, j=1$ follows similarly:
\begin{align*}
\widehat\Psi&\left([\sigma_{2,\alpha,f(X_3)}]^{-1}\,[\sigma_{1,\beta,g(X_2,X_3)}]\,[\sigma_{2,\alpha,f(X_3)}]\right)\\
&=\widehat\Psi\left({[\sigma_{2,\alpha,f(X_3)}]}\right)^{-1}\,\widehat\Psi\left([\sigma_{1,\beta,g(X_2,X_3)}]\right)\,\widehat\Psi\left([\sigma_{2,\alpha,f(X_3)}]\right)\\
&={\mfs_{2,\alpha,f(X_3)}}^{-1}\,\mft_{1,3}\,\mfs_{3,\beta,g(X_2,X_1)}\,\mft_{1,3}\,\mfs_{2,\alpha,f(X_3)}\qquad\text{by (\ref{assignf})}\\
&=\mft_{1,3}\,\mft_{1,3}\,{\mfs_{2,\alpha,f(X_3)}}^{-1}\,\underbrace{\mft_{1,3}\,\mfs_{3,\beta,g(X_2,X_1)}\,\mft_{1,3}}\,\mfs_{2,\alpha,f(X_3)}\,\mft_{1,3}\,\mft_{1,3}\quad\text{since ${\mft_{1,3}}^2=1$}\\
&=\mft_{1,3}\underbrace{\,{\mfs_{2,\alpha,f(X_1)}}^{-1}\,\overbrace{\mfs_{3,\beta,g(X_2,X_1)}}\,\mfs_{2,\alpha,f(X_1)}}\,\mft_{1,3}\qquad\text{by (\ref{special12})}\\
&=\underbrace{\mft_{1,3}\,\overbrace{\mfs_{3,\beta,g(\alpha X_2+f(X_1),X_1)}}\,\mft_{1,3}}\qquad\text{by (\ref{rel2a23})}\\
&=\widehat\Psi\left([\sigma_{1,\beta,g(\alpha X_2+f(X_3),X_3)}]\right)\qquad\text{by (\ref{assignf})}\,,
\end{align*}
completing the proof that the relation (\ref{rel2a}) is respected by $\widehat\Psi$.

Lastly we come to (\ref{rel3a}).  If $\{k,\ell,i\}=\{2,3\}$ then we also have $j\in\{2,3\}$ and we must show that $\mft_{k,\ell}\,\mfs_{i,\alpha,f}\,\mft_{k,\ell}=\mfs_{j,\alpha,f(\tau_{k,\ell})}$.  But this relation holds in $\widetilde\mfH_1$.  Also, if $i=3$ and $\{k,\ell\}=\{1,2\}$, then $j=3$ and the relation holds in $\mfH_2$.  Thus for $i=3$ the only remaining case is $\{k,\ell\}=\{1,3\}$, which follows quickly from (\ref{assignf}).  To wit:
\begin{align*}
\widehat\Psi&([\tau_{1,3}][\sigma_{3,\alpha,f(X_1,X_2)}][\tau_{1,3}])=\widehat\Psi([\tau_{1,3}])\,\widehat\Psi([\sigma_{3,\alpha,f(X_1,X_2)}])\,\widehat\Psi([\tau_{1,3}])\\
&=\mft_{1,3}\,\mfs_{3,\alpha,f(X_1,X_2)}\,\mft_{1,3}\qquad\text{by (\ref{assign23}) and (\ref{Phi[tau]})}\\
&=\widehat\Psi\left([\sigma_{1,\alpha,f(X_3,X_2)}]\right)\qquad\text{by (\ref{assignf})}\\
&=\widehat\Psi\left([\sigma_{1,\alpha,f(\tau_{1,3})}]\right)
\end{align*}
For $i=2$ the remaining cases are $\{k,\ell\}=\{1,2\}$ and $\{k,\ell\}=\{1,3\}$.  For the first:
\begin{align*}
\widehat\Psi&([\tau_{1,2}][\sigma_{2,\alpha,f(X_1,X_3)}][\tau_{1,2}])=\widehat\Psi([\tau_{1,2}])\,\widehat\Psi([\sigma_{2,\alpha,f(X_1,X_3)}])\,\widehat\Psi([\tau_{1,2}])\\
&=\mft_{1,2}\,\underbrace{\mfs_{2,\alpha,f(X_1,X_3)}}\,\mft_{1,2}\qquad\text{by (\ref{assign23}) and (\ref{Phi[tau]})}\\
\intertext{Using the relation $\mfs_{2,\alpha,f(X_1,X_3)}=\mft_{2,3}\,\mfs_{3,\alpha,f(X_1,X_2)}\,\mft_{2,3}$ in $\widetilde\mfH_1\,$:}
&=\lefteqn{\underbrace{\phantom{\mft_{1,2}\,\mft_{2,3}}}}
\mft_{1,2}\,
\lefteqn{\overbrace{\phantom{\mft_{2,3}\,\mfs_{3,\alpha,f(X_1,X_2)}\,\mft_{2,3}}}}
\mft_{2,3}\,\mfs_{3,\alpha,f(X_1,X_2)}\,\underbrace{\mft_{2,3}\,\mft_{1,2}}\\
&=\lefteqn{\overbrace{\phantom{\mft_{1,3}\,\mft_{1,2}}}}
\mft_{1,3}\,\lefteqn{\underbrace{\phantom{\mft_{1,2}\,\mfs_{3,\alpha,f(X_1,X_2)}\,\mft_{1,2}}}}
\mft_{1,2}\,\mfs_{3,\alpha,f(X_1,X_2)}\,\overbrace{\mft_{1,2}\,\mft_{1,3}}\qquad\text{using (\ref{permrel})}\\
\intertext{Using the relation $\mft_{1,2}\,\mfs_{3,\alpha,f(X_1,X_2)}\,\mft_{1,2}=\mfs_{3,\alpha,f(X_2,X_1)}$ in $\mfH_2\,$:}
&=\overbrace{\mft_{1,3}\,\mfs_{3,\alpha,f(X_2,X_1)}\,\mft_{1,3}}\\
&=\widehat\Psi\left([\sigma_{1,\alpha,f(X_2,X_3)}]\right)\qquad\text{by (\ref{assignf})}\\
&=\widehat\Psi\left([\sigma_{1,\alpha,f(\tau_{1,2})}]\right)\,.
\end{align*}
For the second:
\begin{align*}
\widehat\Psi&([\tau_{1,3}][\sigma_{2,\alpha,f(X_1,X_3)}][\tau_{1,3}])=\widehat\Psi([\tau_{1,3}])\,\widehat\Psi([\sigma_{2,\alpha,f(X_1,X_3)}])\,\widehat\Psi([\tau_{1,3}])\\
&=\mft_{1,3}\,\mfs_{2,\alpha,f(X_1,X_3)}\,\mft_{1,3}\qquad\text{by (\ref{assign23}) and (\ref{Phi[tau]})}\\
&=\mfs_{2,\alpha,f(X_3,X_1)}\qquad\text{by (\ref{special12})}\\
&=\widehat\Psi\left([\sigma_{2,\alpha,f(X_3,X_1)}]\right)\\
&=\widehat\Psi\left([\sigma_{2,\alpha,f(\tau_{1,3})}]\right)\,.
\end{align*}
Thus we have verified all the cases when $i=2$ or $i=3$.

Finally we consider $i=1$.  If $\{k,\ell\}=\{2,3\}$, (\ref{rel3a}) is a consequence of (\ref{assignf}):
\begin{align*}
\widehat\Psi&([\tau_{2,3}][\sigma_{1,\alpha,f(X_2,X_3)}][\tau_{2,3}])=\widehat\Psi([\tau_{2,3}])\,\widehat\Psi([\sigma_{1,\alpha,f(X_2,X_3)}])\,\widehat\Psi([\tau_{2,3}])\\
&=\underbrace{\mft_{2,3}\,\mft_{1,3}}\,\mfs_{3,\alpha,f(X_2,X_1)}\,\underbrace{\mft_{1,3}\,\mft_{2,3}}\qquad\text{by (\ref{Phi[tau]}) and  (\ref{assignf})}\\
&=\lefteqn{\overbrace{\phantom{\mft_{1,3}\,\mft_{1,2}}}}
\mft_{1,3}\,
\lefteqn{\underbrace{\phantom{\mft_{1,2}\,\mfs_{3,\alpha,f(X_2,X_1)}\,\mft_{1,2}}}}
\mft_{1,2}\,\mfs_{3,\alpha,f(X_2,X_1)}\,\overbrace{\mft_{1,2}\,\mft_{1,3}}\qquad\text{using (\ref{permrel})}\\
\intertext{Using the relation $\mft_{1,2}\,\mfs_{3,\alpha,f(X_2,X_1)}\,\mft_{1,2}=\mfs_{3,\alpha,f(X_1,X_2)}$ in $\mfH_2\,$:}
&=\mft_{1,3}\,\overbrace{\mfs_{3,\alpha,f(X_1,X_2)}}\,\mft_{1,3}\\
&=\widehat\Psi\left([\sigma_{1,\alpha,f(X_3,X_2)}]\right)\qquad\text{by (\ref{assignf})}\\
&=\widehat\Psi\left([\sigma_{1,\alpha,f(\tau_{2,3})}]\right)\,.
\end{align*}
If $\{k,\ell\}=\{1,3\}$  we have
\begin{align*}
\widehat\Psi&([\tau_{1,3}][\sigma_{1,\alpha,f(X_2,X_3)}][\tau_{1,3}])=\widehat\Psi([\tau_{1,3}])\,\widehat\Psi([\sigma_{1,\alpha,f(X_2,X_3)}])\,\widehat\Psi([\tau_{1,3}])\\
&=\mft_{1,3}\,\mft_{1,3}\,\mfs_{3,\alpha,f(X_2,X_1)}\,\mft_{1,3}\,\mft_{1,3}\qquad\text{by (\ref{Phi[tau]}) and (\ref{assignf})}\\
&=\mfs_{3,\alpha,f(X_2,X_1)}\qquad\text{since ${\mft_{1,3}}^2=1$}\\
&=\widehat\Psi\left([\sigma_{3,\alpha,f(X_2,X_1)}]\right)\qquad\text{by (\ref{assignf})}\\
&=\widehat\Psi\left([\sigma_{3,\alpha,f(\tau_{1,3})}]\right)\,.
\end{align*}
If $\{k,\ell\}=\{1,2\}$  we have
\begin{align*}
\widehat\Psi&([\tau_{1,2}][\sigma_{1,\alpha,f(X_2,X_3)}][\tau_{1,2}])=\widehat\Psi([\tau_{1,2}])\,\widehat\Psi([\sigma_{1,\alpha,f(X_2,X_3)}])\,\widehat\Psi([\tau_{1,2}])\\
&=\underbrace{\mft_{1,2}\,\mft_{1,3}}\,\mfs_{3,\alpha,f(X_2,X_1)}\,\underbrace{\mft_{1,3}\,\mft_{1,2}}\qquad\text{by (\ref{Phi[tau]}) and (\ref{assignf})}\\
&=\lefteqn{\overbrace{\phantom{\mft_{1,3}\,\mft_{2,3}}}}
\mft_{1,3}\,
\lefteqn{\underbrace{\phantom{\mft_{2,3}\,\mfs_{3,\alpha,f(X_2,X_1)}\,\mft_{2,3}}}}
\mft_{2,3}\,\mfs_{3,\alpha,f(X_2,X_1)}\,\overbrace{\mft_{2,3}\,\mft_{1,3}}\qquad\text{using (\ref{permrel})}\\
\intertext{Using the relation $\mft_{2,3}\,\mfs_{3,\alpha,f(X_2,X_1)}\,\mft_{2,3}=\mfs_{2,\alpha,f(X_3,X_1)}$ in $\widetilde\mfH_1\,$:}
&=\mft_{1,3}\,\overbrace{\mfs_{2,\alpha,f(X_3,X_1)}}\,\mft_{1,3}\\
&=\mfs_{2,\alpha,f(X_1,X_3)}\qquad\text{by (\ref{special12})}\\
&=\widehat\Psi\left([\sigma_{2,\alpha,f(X_1,X_3)}]\right)\\
&=\widehat\Psi\left([\sigma_{2,\alpha,f(\tau_{1,2})}]\right)\,.
\end{align*}
The proof of Theorem \ref{main} is now complete.

\section{Concluding Remarks and Questions} 

The combinatoric upshot of Theorem \ref{main} is that $\TA_3(K)$ is the colimit of a ``triangle of groups" in Stallings' sense (see \cite{Stall}), comprising $\widetilde H_1$, $H_2$, and $H_3$, their pairwise intersections, and the intersection of all three.  These groups form the stabilizers of the three vertices, three edges and face, respectively, of a simplex $\mathfrak f$ in a simply connected complex of triangles $\mathcal D$ on which $\TA_3(K)$ acts, and for which $\mathfrak f$ serves as a fundamental domain.  There are unanswered questions about $\mathcal D$.  For example, is it $2$-connected (i.e., does every continuous image of the 2-sphere in $\mathcal D$ contract in $\mathcal D$), and does it have infinite diameter (i.e., is the number of faces need to connect two arbitrary points unbounded)?

%The $3$-dimensional simplicial complex $\mathcal D$ can be realized as follows:  More generally we construct an $n$-dimensional simplicial complex $\mathcal E_n$ whose ertices are rank $(i+1)$ vector spaces $V$ in $K[\X]$ containing $K$, where $1\le i\le n$ such that $k[\X]=k[V]^{[n-i]}$.  Vertices $V_1,\ldots,V_r$ of strictly ascending rank form an $r$-simplex if $V_1\subset\cdots\subset V_r$.  There is an obvious action of $\GA_n(K)$ on $\mathcal E_n$.  Note that $\mathcal E_n$ contains the $n$-simplex $\mathfrak s$ determined by the $n$ vertices $V_i$, $1\le i\le n$, defined by (\ref{vecsp}) in \S\ref{S4}.  This is a fundamental domain for the action, and the subgroup $H_i$ is the stabilizer of $V_i$, by its definition (\ref{Hdef}).  For $n=2$ this is the tree which gives the structure Theorem for $\GA_2(K)$ (Example \ref{JVdK}).  For $n\ge3$ we do not know if $\mathcal E_n$ is connected, or simply connected.
%
%What if we remove the tilde from $H_1$?
%
%Contain spheres?  (mention the $\Cr_2$ situation, Lami's result)
%
%
%
%Infinite diameter? 
%
%Tree of $H_1$ lies in the complex.

\bibliographystyle{amsplain}
\bibliography{Refs}

\vskip 15pt

\noindent{\small \sc Department of Mathematics, Washington University in St.
Louis, St. Louis, MO 63130 } {\em E-mail}: wright@math.wustl.edu

\end{document}